\documentclass[12pt]{article}
\begin{document}
\title{An Analysis of the Min-max Algorithm}
\author{Jerzy Cis{\l}o \\
Institute of Theoretical Physics  \\
University of Wroc{\l}aw  \\
pl. M.Borna 9, 50-205 Wroc{\l}aw, Poland  \\
cislo@ift.uni.wroc.pl
}
\date{}

\maketitle

\begin{abstract}
We consider the matrix $A_{ij}$, whose elements are independent random variables.
We calculate the mean value of the number of the elements that we need to read to find
$\min_{i}\max_{j} A_{ij}$.
\end{abstract}

\section{The Min-max Algorithm}
In this paper we are going to analyze the process of finding $\min_i\max_j A_{ij}$ [1,2,3], 
for which the reading of any single matrix element is particularly time consuming. 
However, this process--used in game theory and economy--may be shortened as 
to find the value of the Min-max, we do not have to read the whole matrix.
Let us look at the following example:
 $$
\matrix{   5  &  3  &  {\bf 14}  &  8  \cr
           2  &  4  &   {\bf 9}  &  1  \cr
           7  & {\bf 16}  &  \underline{13}  & \underline{11}  \cr
          12  &  \underline{6}  &  \underline{10}  & \underline{\bf 15}  }
$$
The maximal elements in each row are printed bold, and the elements which we do not have to read are
underlined.

\newpage

The algorithm of finding Min-max of matrix $n \times k$ reads
 \begin{verbatim}

minmax= infinity
FOR i=1 TO n
    max= -infinity
    FOR j=1 TO k
        t=A[i][j]
        IF t > max THEN max=t
        IF t > minmax THEN BREAK
    IF max < minmax THEN minmax=max

\end{verbatim}
The essential instruction in this algorithm is: 
   "IF t $>$ minmax THEN BREAK".
 
The presented algorithm is a restricted version of algorithm alpha-beta pruning,
that is used in logical computer games such as chess, Othello, etc.
 
The search function of algorithm alpha-beta pruning reads [4,5]

 \begin{verbatim}
SEARCH(n,alpha)
    IF n = 0 THEN RETURN  - evaluation
    beta = - infinity
    FOR j=1 TO m  (list of moves)
        MOVE(j)
        t = - SEARCH( n-1, - beta)
        UNDO MOVE
        IF t <= alpha THEN RETURN beta
        IF t < beta THEN beta=t
    RETURN beta
\end{verbatim}

\section{Two equivalent models}
Generally, the time of performance of the algorithm depends on the distribution of matrix elements.

In this paper we consider the probability model in which we assume that the matrix elements $(A_{11}, ..., A_{nk})$ 
are independent uniformly distributed random variables on $\langle 0,1 \rangle $.

Our aim is to find the mean value of the number of matrix elements that are read while searching for the min-max.

We can also consider the discrete model in which the matrix elements
$(A_{11}, ..., A_{nk})$ are permutations of the set $(1,2,..., nk)$, each permutation being of the same probability.

For both models, the mean values are the same.

Still another model is the discrete model in which the matrix elements belong
to the finite set (i.e. the elements can be repeated). For this model, however, the obtained mean value 
is different from the ones achieved for the first two models.

\section{Recursion}

When the algorithm begins to operate, the parameter $a$ takes the maximum value.
After each row the parameter $a$ can decrease.
The further search depends on the value of $a$.
Let $M(a)$ stand for the mean value of the read elements, on condition that
the algorithm begins to operate with the parameter $a \in \langle 0,1 \rangle $.

Now at the beginning of the matrix
$A_{ij}$, let us add a new row $( x_{1}, x_{2}, ...x_{k})$ consisting of random independent variables uniformly distributed
on $\langle 0,1 \rangle $.

Let us take an arbitrary number $a \in \langle 0,1 \rangle$.

The probability that  $x_1,x_2, x_{j-1} \leq a < x_j$ is equal to $(1-a)^{j-1}$.

If $x_1,x_2, \dots, x_{j-1} \leq a < x_j$, the algorithm reads the elements
 $x_1,\dots ,x_j$, and, on average, $M(a)$ elements from the remaining part of the matrix. 

The probability that $\{x_1,\dots,x_k\}$ are in $\langle t, t+\Delta\rangle $ is equal 
to $(1+\Delta)^k-t^k$.

If $\{x_1,\dots,x_k\} \subseteq \langle t, t+\Delta \rangle$, the algorithm reads all $k$
elements from the added row, and, on average, $M(t)$ elements from the remaining part of the matrix.

Therefore the mean number of the read matrix elements is equal to

$$ \tilde{M}(\alpha)= (1-\alpha)(1+M(\alpha)) +
                   (1-\alpha)(2+M(\alpha)) +...
                  +(1-\alpha)(k+M(\alpha)) + $$
$$                +  \int_{0}^{\alpha} (k+M(t))
                  {d t^{k} \over dt} \; dt 
 = 1+ \alpha + ... + \alpha^{k-1} +M(\alpha)+
   \int_{0}^{\alpha} t^{k} {d M(t) \over dt} \;dt. $$

\section{Calculations}
The obtained recursion makes it possible to calculate the mean number of the read matrix elements. 

Let $M_{nk}(a)$ stand for the mean number of the read  matrix elements, on condition that
the algorithm begins to operate with parameter  $a \in \langle 0,1\rangle$,
$n$ and $k$ representing, respectively, the number of columns and number of rows in the matrix.

Using the relation from the previous section we can write 
\begin{equation}\label{rec}
M_{nk}(a)=M_{n-1,k}(a)-\int_{0}^{a}t^{k}M_{n-1,k}(t) \;dt
\end{equation}
$$ + 1+a+ a^{2}+...+a^{k-1}. $$
This relation, together with  the initial condition $M_{0k}(a)=0$,
allows us to find  $M_{nk}(a)$ for any $n$.

Differentiating the equation (\ref{rec}) we obtain
\begin{equation} \label{rel}
\frac{dM_{nk}(a)}{da} = (1-a^k) \frac{dM_{n-1,k}(a)}{da} + (1+2a+3a^2+\cdots +(k-1)a^{k-1}).
\end{equation}
We remember that
\begin{equation}
\frac{dM_{0k}(a)}{da}=0, \quad M_{nk}(0)=0.
\end{equation}
Iterating $n$ times the relation (\ref{rel})  we get
\begin{equation}
\frac{dM_{nk}(a)}{da}=\sum_{j=0}^{n-1} (1-a^k)^j (1+2a+3a^2+\cdots +(k-1)a^{k-1})
\end{equation}
$$
=\frac{1-(1-a^k)^n}{a^k}(1+2a+3a^2+\cdots +(k-1)a^{k-1}).
$$
Hence
\begin{equation}\label{resa}
M_{nk}(a)= \int_0^a \frac{1-(1-t^k)^n}{t^k}(1+2t+3t^2+\cdots +(k-1)t^{k-1}) \; dt=
\end{equation} 
$$
M_{nk}(a)= n \sum_{l=o}^{k-1} a^{k} +
\sum_{l=1}^{k-1} \sum_{j=1}^{n-1}
\pmatrix{ n \cr j+1} (-1)^{j} {l \over l+jk} a^{l+jk}.
$$
\section{Results}
We obtain our main result, that is the formula for the mean number of the read
matrix elements,  by substituting $a=1$ in the formula (\ref{resa})
\begin{equation}
M_{nk}(1)= n k + \sum_{l=1}^{k-1} \sum_{j=1}^{n-1}
\pmatrix{ n \cr j+1} (-1)^{j} {l \over l+jk}.
\end{equation}
Using the formula [6]
\begin{equation}
F_{n}(x)=n+ \sum_{j=1}^{n-1} \pmatrix { n \cr j+1}
(-1)^{j} {x \over j+x} = 
\end{equation}
$$ 
1 + {1 \over (1+x)} + {2! \over (1+x)(2+x) } +
\cdots + {(n-1)! \over (1+x)(2+x)...(n-1 +x)},
$$
we can rewrite our result in the following form:
\begin{equation}
M_{nk}(1)=\sum_{l=0}^{k-1} F_{n}(l / k) =
   n+k-1+ \sum_{l=1}^{k-1} \sum_{j=1}^{n-1} \prod_{r=1}^{j}
   {rk \over rk + l}.
\end{equation}
For small values of $n$ and $m$, we have
$$ M_{22}(1)=3{2 \over 3}, \quad M_{23}(1)=5{7 \over 20}, \quad M_{24}(1)=7{4 \over 105}, $$
$$ M_{32}(1)=5{1 \over 5}, \quad M_{33}(1)=7{31 \over 70}, \quad   M_{34}(1)=9{2419 \over 3465}, $$ 
$$ M_{42}(1)=6{23 \over 35}, \quad  M_{43}(1)=9{30 \over 77}, \quad M_{44}(1)=12{6493 \over 45045}. $$
Also, we would like to add that the following inequalities may be checked:
\begin{equation}
{M_{nk}(1) \over n} > {M_{n+1,k}(1) \over n+1}, \quad
{M_{nk}(1) \over k} > {M_{n,k+1}(1) \over k+1},
\end{equation}
There is a good approximation for $M_{nk}(1)$
$$ M_{nk}(1) \approx  \frac{n}{2}+k\left(\frac{n}{\log n}+ \frac{n}{\log^2 n} \right). $$  
The error of the formula is less then 3 percent for $n \geq 28 ,k \geq 17$ ,
and less then 1 percent for $n \geq 38, k \geq 20$.

\end{document}